\documentclass{IEEEtran}
\usepackage{graphicx}
\usepackage{amsmath}
\usepackage[exponent-product = \cdot]{siunitx}
\usepackage{tikz}
\usepackage{circuitikz}
\usepackage{pgfplots}
\usepgfplotslibrary{groupplots}
\usepackage[figuresright]{rotating}
\usepackage{color,xcolor}
\usepackage{amsfonts,amssymb,mathtools}
\usepackage{amstext}
\usepackage{upgreek}
\usepackage[tight,footnotesize]{subfigure}
\usepackage[tight,nice]{units}
\usepackage[figure]{algorithm2e}
\usepackage[final]{changes}
\usepackage[mathcal]{eucal}
\usepackage{pgfplots}
\usepackage{verbatim}
\usepackage{xspace}
\usepackage[square,numbers]{natbib}
\catcode`@=11
\def\frownfill{$\scriptscriptstyle\m@th\mathord\frown$}
\def\bow#1{\vbox{\m@th\ialign{##\crcr
      \hfil\frownfill\hfil\crcr\noalign{\kern-0.2\p@\nointerlineskip}
      $\hfil\displaystyle{#1}\hfil$\crcr}}}
\def\bbow#1{\vbox{\m@th\ialign{##\crcr
     \hfil\frownfill\hfil\crcr\noalign{\kern-0.7\p@\nointerlineskip}
     \hfil\frownfill\hfil\crcr\noalign{\kern-0.3\p@\nointerlineskip}
      $\hfil\displaystyle{#1}\hfil$\crcr}}}
\def\widefrownfill{$\m@th\mathord\frown$}
\def\widebow#1{\vbox{\m@th\ialign{##\crcr
      \hfil\widefrownfill\hfil\crcr\noalign{\kern-0.9\p@\nointerlineskip}
      $\hfil\displaystyle{#1}\hfil$\crcr}}}
\def\widebbow#1{\vbox{\m@th\ialign{##\crcr
     \hfil\widefrownfill\hfil\crcr\noalign{\kern-1.8\p@\nointerlineskip}
     \hfil\widefrownfill\hfil\crcr\noalign{\kern-0.9\p@\nointerlineskip}
      $\hfil\displaystyle{#1}\hfil$\crcr}}}

\newcommand{\Bfield}{\ensuremath{\vec{B}}} 
\newcommand{\Dfield}{\ensuremath{\vec{D}}} 
\newcommand{\Efield}{\ensuremath{\vec{E}}} 
\newcommand{\Hfield}{\ensuremath{\vec{H}}} 
\newcommand{\Jfield}{\ensuremath{\vec{J}}}

\newcommand{\eps}{\ensuremath{\varepsilon}}

\newcommand{\ve}{\protect\bow{\mathrm{\mathbf{e}}}}

\newcommand{\vh}{\protect\bow{\mathrm{\mathbf{h}}}}

\newcommand{\fb}{\protect\bbow{\mathrm{\mathbf{b}}}}

\newcommand{\fd}{\protect\bbow{\mathrm{\mathbf{d}}}}

\newcommand{\fj}{\protect\bbow{\mathrm{\mathbf{j}}}}

\newcommand{\qfit}{\mathrm{\mathbf{q}}}

\newcommand{\Meps}{\mathbf{M}_{\varepsilon}}

\newcommand{\Msigma}{\mathbf{M}_{\sigma}}
\newcommand{\Mmu}{\mathbf{M}_{\mu}}

\newcommand{\dw}{\Delta w}

\newcommand{\divfit}{\mathbf{S}}
\newcommand{\divdfit}{\widetilde{\mathbf{S}}}
\newcommand{\curlfit}{\mathbf{C}}
\newcommand{\curldfit}{\widetilde{\mathbf{C}}}

\newcommand{\T}{\mathbf{T}}

\newcommand{\real}{\ensuremath{\mathbb{R}}}

\newcommand{\complex}{\ensuremath{\mathbb{C}}}

\newcommand{\curl}{\nabla \times}
\renewcommand{\div}{\nabla \cdot}

\renewcommand{\vec}[1]{\mathbf{#1}}

\newcommand{\bu}{\mbox{\boldmath$u$}}
\newcommand{\bv}{\mbox{\boldmath$v$}}
\newcommand{\bw}{\mbox{\boldmath$w$}}

\newcommand{\bb}{\mbox{\boldmath$b$}}

\newcommand{\magenergy}[2]{\bigl\langle\vh^{(#1)},\vh^{(#2)}\bigr\rangle_\mu}
\newcommand{\elecenergy}[2]{\bigl\langle\ve^{(#1)},\ve^{(#2)}\bigr\rangle_\varepsilon}
\newcommand{\tr}[1]{\bigl({#1}\bigr)^{\!\!\top}\!}
\newcommand{\MM}{\mathbf{M}}
\newcommand{\KK}{\mathbf{K}}
\newcommand{\A}{\mathbf{A}}

\newcommand{\cstmws}{\textsc{CST MICROWAVE STUDIO\textsuperscript{\textregistered}} \xspace}

\definecolor{tud1a}{RGB}{93,133,195}
\definecolor{tud2a}{RGB}{0,156,218}
\definecolor{tud3a}{RGB}{80,182,149}
\definecolor{tud4a}{RGB}{175,204,80}
\definecolor{tud5a}{RGB}{221,223,72}
\definecolor{tud6a}{RGB}{255,224,92}
\definecolor{tud7a}{RGB}{248,186,60}
\definecolor{tud8a}{RGB}{238,122,52}
\definecolor{tud9a}{RGB}{233,80,62}
\definecolor{tud10a}{RGB}{201,48,142}
\definecolor{tud11a}{RGB}{128,69,151}

\definecolor{tud1b}{RGB}{0,90,169}
\definecolor{tud2b}{RGB}{0,131,204}
\definecolor{tud3b}{RGB}{0,157,129}
\definecolor{tud4b}{RGB}{153,192,0}
\definecolor{tud5b}{RGB}{201,212,0}
\definecolor{tud6b}{RGB}{253,202,0}
\definecolor{tud7b}{RGB}{245,163,0}
\definecolor{tud8b}{RGB}{236,101,0}
\definecolor{tud9b}{RGB}{230,0,26}
\definecolor{tud10b}{RGB}{166,0,132}
\definecolor{tud11b}{RGB}{114,16,133}

\definecolor{tud1c}{RGB}{0,78,138}
\definecolor{tud2c}{RGB}{0,104,157}
\definecolor{tud3c}{RGB}{0,136,119}
\definecolor{tud4c}{RGB}{127,171,22}
\definecolor{tud5c}{RGB}{177,189,0}
\definecolor{tud6c}{RGB}{215,172,0}
\definecolor{tud7c}{RGB}{210,135,0}
\definecolor{tud8c}{RGB}{204,76,3}
\definecolor{tud9c}{RGB}{185,15,34}
\definecolor{tud10c}{RGB}{149,17,105}
\definecolor{tud11c}{RGB}{97,28,115}

\definecolor{tud1d}{RGB}{36,53,114}
\definecolor{tud2d}{RGB}{0,78,115}
\definecolor{tud3d}{RGB}{0,113,94}
\definecolor{tud4d}{RGB}{106,139,55}
\definecolor{tud5d}{RGB}{153,166,4}
\definecolor{tud6d}{RGB}{174,142,0}
\definecolor{tud7d}{RGB}{190,111,0}
\definecolor{tud8d}{RGB}{169,73,19}
\definecolor{tud9d}{RGB}{156,28,38}
\definecolor{tud10d}{RGB}{115,32,84}
\definecolor{tud11d}{RGB}{76,34,106}

\pgfplotscreateplotcyclelist{temfbarplot cycle}{
	{draw=tud1d,		fill=tud1b},
	{draw=tud3d,		fill=tud3b},
	{draw=tud5d,		fill=tud5b},
	{draw=tud7d,		fill=tud7b},
	{draw=tud9d,		fill=tud9b},
	{draw=tud11d,	 fill=tud11b}, 
	{draw=tud2d,		fill=tud2b},
	{draw=tud4d,		fill=tud4b},
	{draw=tud6d,		fill=tud6b},
	{draw=tud8d,		fill=tud8b},
	{draw=tud10d,	 fill=tud10b},
}
\pgfplotscreateplotcyclelist{temflineplot cycle}{
	{tud1b, mark=*, mark size=1.5pt, thick},
	{tud3b, mark=square*, mark size=1.3pt, thick},
	{tud5b, mark=triangle*, mark size=1.5pt, thick},
	{tud7b, mark=diamond*, mark size=1.5pt, thick},
	{tud9b, mark=*, mark size=1.5pt, thick},
	{tud11b, mark=square*, mark size=1.3pt, thick},
	{tud2b, mark=triangle*, mark size=1.5pt, thick},
	{tud4b, mark=diamond*, mark size=1.5pt, thick},
	{tud6b, mark=*, mark size=1.5pt, thick},
	{tud8b, mark=square*, mark size=1.3pt, thick},
	{tud10b, mark=triangle*, mark size=1.5pt, thick},
}

\pgfplotsset{
	colormap={temfcolormap}{
		color(0cm)=(tud1b);
		color(1cm)=(tud6b);
		color(2cm)=(tud8b);
		color(3cm)=(tud9b)
	},
	temf3dplot/.style={
		temfbaseplot,
		xmajorgrids=true,
		ymajorgrids=true,
		major grid style={dotted},
		axis line style = thin,
		axis x line*=box,
		axis y line*=box,
	},
	temflineplot/.style={
		temfbaseplot,
		xmajorgrids=true,
		ymajorgrids=true,
		major grid style={dotted},
		axis line style = thin,
		axis x line*=bottom,
		axis y line*=left,
		legend style={cells={anchor=west}},
		cycle list name=temflineplot cycle,
	},
	temfbarplot base/.style={
		temfbaseplot,
		bar width=6pt,
		axis y line*=none,
	},
	temfbarplot/.style={
		temfbarplot base,
		ybar,
		xmajorgrids=false,
		ymajorgrids=true,
		area legend,
		legend image code/.code={\draw[#1] (0cm,-0.1cm) rectangle (0.15cm,0.1cm);},
		cycle list name=temfbarplot cycle,
	},
	horizontal temfbarplot/.style={
		temfbarplot base,
		xmajorgrids=true,
		ymajorgrids=false,
		xbar stacked,
		area legend,
		legend image code/.code={\draw[#1] (0cm,-0.1cm) rectangle (0.15cm,0.1cm);},
		cycle list name=temfbarplot cycle,
	},
	temfbaseplot/.style={
		colormap name=temfcolormap,
		legend style={fill=white,cells={anchor=west}},
		x tick label style={font=\footnotesize},
		y tick label style={font=\footnotesize},
		legend style={font=\footnotesize},
		major grid style={dotted},axis x line*=bottom
	},
	disable thousands separator/.style={/pgf/number format/.cd,1000 sep={}},
}

\usepgfplotslibrary{groupplots}
\usetikzlibrary{decorations.markings}
\usetikzlibrary{arrows}
\usetikzlibrary{patterns}
\pgfplotsset{cycle list={
tud9b,solid\\
tud9b,dashed\\
tud5b,dashdotted\\
black,dotted\\
brown,dashdotdotted\\
teal,solid\\
orange,dashed\\
violet,dotted\\
cyan,dashdotted\\
tud5d,dashdotdotted\\
magenta,solid\\
gray,dashed\\}
}

\begin{document}
\title{ParaExp using Leapfrog as Integrator for High-Frequency Electromagnetic Simulations}
\author{\IEEEauthorblockN{%
M. Merkel\IEEEauthorrefmark{1}, %
I. Niyonzima\IEEEauthorrefmark{1}$^{,}$\IEEEauthorrefmark{2}, and %
S. Sch{\"o}ps\IEEEauthorrefmark{1}$^{,}$\IEEEauthorrefmark{2}}%

\IEEEauthorblockA{\IEEEauthorrefmark{1}Graduate School of Computational Engineering (GSC CE), Technische Universit\"{a}t Darmstadt, Germany.}

\IEEEauthorblockA{\IEEEauthorrefmark{2}Institut f\"{u}r Theorie Elektromagnetischer Felder (TEMF), Technische Universit\"{a}t Darmstadt, Germany.}
\thanks{Corresponding author: Melina Merkel,Technische Universit{\"a}t Darmstadt
Institut f{\"u}r Theorie Elektromagnetischer Felder, Schlo{\ss}gartenstr. 8
64289 Darmstadt (anna\_melina.merkel@stud.tu-darmstadt.de)}%
}
\maketitle
\begin{abstract}
        Recently, ParaExp was proposed for the time integration of linear hyperbolic problems. It splits the time interval of interest into
        sub-intervals and computes the solution on each sub-interval in parallel. The overall solution is decomposed into a particular
        solution defined on each sub-interval with zero initial conditions and a homogeneous solution propagated by the matrix exponential
        applied to the initial conditions. The efficiency of the method depends on fast approximations of this matrix exponential based on recent results
        from numerical linear algebra. This paper deals with the application of ParaExp in combination with Leapfrog to electromagnetic wave problems in time-domain. Numerical
        tests are carried out for a simple toy problem and a realistic spiral inductor model discretized by the Finite Integration Technique.
\end{abstract}


{\color{black}
    
\section{Introduction}

The simulation of high-frequency electromagnetic problems is often carried out in 
frequency domain. This choice is motivated by the linearity of the underlying 
governing equations. However, the solution of problems in frequency domain may 
require the resolution of very large linear systems of equations and this becomes 
particularly inconvenient for broadband simulations such that approximations like 
model order reduction are typically used, e.g. \cite{slone03,Floch_2015aa,Paquay_2016aa}. 
The coupling with nonlinear time-dependent systems and the computation of transients 
are other cases where time-domain simulations outperform frequency-domain simulations.

On the other hand, the numerical complexity resulting from time-domain simulations 
may also become prohibitively expensive. Parallelization in `space', e.g., matrix-vector 
multiplications corresponding to the application of the curl operator, using multicore architectures is well established in academic and industrial software environments \cite{CST_2016aa}. However, the parallelization efficiency eventually saturates with increasing number of cores depending on the memory bandwidth of the involved hardware. Time-domain parallelization is a promising extension to domain decomposition in space.

The development and application of parallel-in-time methods dates back more 
than 50 years, see \cite{nievergelt1964parallel}. These methods can be direct 
\cite{christlieb2010parallel,gander2013paraexp} or iterative 
\cite{lions2001parareal,minion2011hybrid}. They can also be well suited for small 
scale parallelization \cite{miranker1967parallel,womble1990time} or large parallelization 
\cite{gander2013paraexp,minion2011hybrid}. Recently, the \emph{Parareal method} 
gained interest \cite{lions2001parareal,gander2007analysis}. In its initial version, 
Parareal was developed for large scale semi-discretized parabolic partial differential 
equations (PDEs). It involves the splitting of the time interval and the resolution 
of the governing ordinary differential equation (ODE) in parallel on each sub-interval 
using a fine propagator which can be any classical time-stepper with a fine time grid.
A coarse propagator distributes the initial conditions for each sub-interval during 
the Parareal iterations. It is typically obtained by a time stepper with a coarse 
grid on the entire time interval. Parareal iterates the resolution of both the 
coarse and the fine problems until convergence.

Most parallel-in-time methods fail for hyperbolic problems. In the case of Parareal, 
analysis has shown that it may lead to the \emph{beating} phenomenon depending on 
the structure of the system matrix \cite{farhat2006beat}. It may even become unstable 
if the eigenvalues of the matrix are purely imaginary which is the case in the 
presence of undamped electromagnetic waves.

In this paper we apply the \emph{ParaExp method} from \cite{gander2013paraexp} for 
the parallelization of time-domain resolutions of hyperbolic equations that govern 
the electromagnetic wave problems as initially proposed in \cite{merkel2016}.

The method splits the time interval into sub-intervals 
and solves smaller problems on each sub-interval as visualized in Figure \ref{fig:timedecomposition}. 
Using the theory of linear ordinary differential equations, the total solution for 
each sub-interval is decomposed into particular solution with zero initial conditions 
and homogeneous solutions with initial conditions from previous intervals.
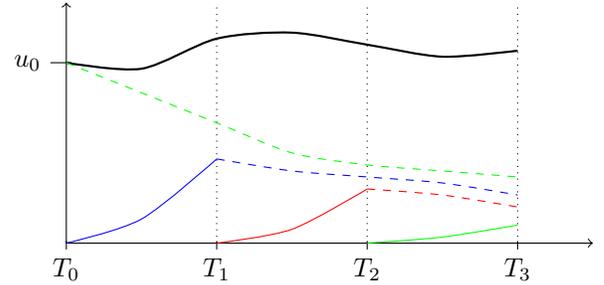
\begin{figure}[t!]
\centering
    \begin{tikzpicture} [yscale=0.4,scale=2]
\draw [->] (0,0)--(3.5,0);
\draw [->] (0,0)--(0,4);
\foreach \T in {0,...,3} \draw (\T,1pt) -- (\T,-3pt) node[anchor=north] {$T_\T$};
\foreach \T in {1,...,3} \draw [dotted](\T,0) -- (\T,4);
\draw (1pt,3) -- (-3pt,3) node[anchor=east] {$u_0$};

\draw [thick] plot [smooth] coordinates{(0,3)(0.5,2.9)(1,3.4)(1.5,3.5)(2,3.3)(2.5,3.1)(3,3.2)};
\draw [blue] plot [smooth] coordinates {(0,0)(0.5,0.4)(1,1.4)};
\draw [red] plot [smooth] coordinates {(1,0)(1.5,0.23)(2,0.9)};
\draw [green] plot [smooth] coordinates {(2,0)(2.5,0.1)(3,0.3)};

\draw [dashed,blue] plot [smooth] coordinates {(1,1.4)(1.5,1.2)(2,1.1)(2.5,1.0)(3,0.8)};
\draw [dashed,red] plot [smooth] coordinates {(2,0.9)(2.5,0.8)(3,0.6)};
\draw [dashed,green] plot [smooth] coordinates {(0,3)(0.5,2.5)(1,2)(1.5,1.5)(2,1.3)(2.5,1.2)(3,1.1)};
\end{tikzpicture}
     \vspace{-0.2cm}
    \caption{Schematic view of the decomposition of time and solution. Vertical 
    dotted lines denote the sub-intervals, solid lines represent the solution of 
    the inhomogeneous sub-problems and dashed lines represent the solution of 
    the homogeneous sub-problems. {The thick black line represents the overall solution.} Colors indicate the employed processors, cf. \cite{gander2013paraexp}}
    \label{fig:timedecomposition}
    \vspace{-0.2cm}
\end{figure}

The paper is organized as follows: in Section \ref{sec:maxwell-equations} we introduce 
Maxwell's equations and derive the governing system of ODEs for the wave equation 
obtained by the Finite Integration Technique (FIT). This system is then used in 
Section \ref{sec:paraexp} for the presentation of the ParaExp method following the 
lines of \cite{gander2013paraexp}. The mathematical framework is briefly sketched 
and the details of the algorithm are discussed. The combination of ParaExp with 
Leapfrog is proposed. Section \ref{sec:applications} deals with numerical examples. 
We consider two applications: a simple wave guide problem and a realistic spiral 
inductor model discretized by the Finite Integration Technique. The examples are 
investigated in terms of efficiency, energy conservation and frequency spectrum.     \section{Space and Time Discretization of Maxwell's equations}
\label{sec:maxwell-equations}

In an open, bounded domain $\Omega \subset \mathbb{R}^3$ and 
$t \in \mathcal{I} = (t_0,t_\mathrm{end}] \subset \mathbb{R}_{\geq0}$, the evolution 
of electromagnetic fields is governed by Maxwell's equations on $\Omega \times \mathcal{I}$,
see e.g. \cite{jackson1999classical}:
\begin{align}
	\curl{\Efield}  &= -\partial_t \Bfield ,
	&\quad
	\curl{\Hfield}  &= \partial_t \Dfield+\Jfield,
	\\
	\div{\Bfield} &= 0,
	&\quad
	\div{\Dfield} &= \rho
	\label{eq:maxwell-equations}
\end{align}
with suitable initial and boundary conditions at time $t_0$ and $\partial\Omega$, respectively.
In presence of linear materials, these equations are completed by constitutive 
laws \cite{jackson1999classical}:
\begin{align}
    \Dfield=\varepsilon\Efield, \quad 
    \Jfield=\sigma \Efield + \Jfield_\text{s}, \quad 
    \Bfield=\mu\Hfield.
  \label{eq:const-laws-1}
\end{align}
In these equations, $\Hfield$ is the magnetic field [A/m], $\Bfield$ the magnetic 
flux density [T], $\Efield$ the electric field [V/m], $\Dfield$ the electric flux 
density [C/m$^2$], $\Jfield$, $\Jfield_{\sigma}=\sigma\Efield$, $\Jfield_{\rm d}=\partial_t \Dfield$ and $\Jfield_s$ 
are the total, Ohmic, displacement and electric source current densities [A/m$^2$], $\rho$ 
is the electric charge density [C/m$^3$]. The material properties $\sigma$, $\varepsilon$ 
and $\mu$ are the electric conductivity, the electric permittivity and the magnetic 
permeability, respectively. In this paper, we consider electromagnetic wave 
propagation in non-conducting media which are free of charges, i.e., $\sigma \equiv 0$ 
and $\rho \equiv 0$.

The space discretization of Maxwell's equations \eqref{eq:maxwell-equations}-\eqref{eq:const-laws-1} 
using the Finite Integration Technique (FIT) \cite{Weiland_1977aa, Weiland_1996aa} 
on a staggered grid pair with primal $n=n_{x}\cdot n_{y} \cdot n_{z}$ grid points leads to the equations
\begin{align}
	\label{eq:maxwell_grid1}
	\curlfit \ve &= - \operatorname{d}_t\fb,
	&\quad
	\curldfit \vh &= \operatorname{d}_t\fd+\fj,
	\\
	\label{eq:maxwell_grid2}
	\divdfit\fb &=\mathbf{0},
	&\quad
	\divfit\fd&=\qfit
\end{align}
where $\curlfit$, $\curldfit\in\real^{n_\text{dof}\times n_\text{dof}}$ are the 
discrete curl operators, $\divfit$, $\divdfit\in\real^{n\times n_\text{dof}}$ the 
discrete divergence operators, which are all defined on the primal and dual grid, 
respectively ($n_\text{dof}\approx 3n$). The fields are semi-discretely given by 
$\ve$, $\vh$, $\fd$, $\fj$, $\fb:\mathcal{I}\to\real^{n_\text{dof}}$ and 
$\qfit:\mathcal{I}\to\real^{n}$ corresponding to electric and magnetic voltages, 
electric fluxes, electric currents, magnetic fluxes and electric charges, respectively. 
They are linked by the material relations
\begin{align}
	\label{eq:maxwell_grid3}
	\fd=\Meps\ve, \quad
	\fj=\Msigma\ve+\fj_\text{s}, \quad
	\fb=\Mmu\vh
\end{align}
where $\Meps$ and $\Mmu$ are diagonal positive-definite material matrices of 
permittivities $\eps$ and permeabilities $\mu$. The conductivity 
matrix $\Msigma$ will not be considered as mentioned above.

The system (\ref{eq:maxwell_grid1}-\ref{eq:maxwell_grid3}) can be rewritten as an 
initial value problem (IVP)
\begin{align}
    \MM \operatorname{d}_t \bar{\mathbf{u}} + \KK \bar{\mathbf{u}} &= \bar{\mathbf{g}}(t), & \bar{\mathbf{u}}(t_0)&= \bar{\mathbf{u}}_0.
\label{eq:FIT_ode_1}
\end{align}
with unknown voltages $\bar{\mathbf{u}}^{\top} := [\vh^{\top},\ve^{\top}]$, given excitation $\bar{\mathbf{g}}^{\top}:=[\mathbf{0},\fj_\text{s}^{\top}]$ and the matrices
\begin{align}
    \MM := 
    \begin{bmatrix}
        \Mmu &  \boldsymbol{0} \\
        \boldsymbol{0}               & \Meps  
    \end{bmatrix}
    \, , \quad 
    \KK := 
    \begin{bmatrix}
        \boldsymbol{0}  & \curlfit \\
        -\curldfit       & \boldsymbol{0}  
    \end{bmatrix}. \label{eq:fit_matrices}
\end{align}
Exploiting a similarity transformation by the matrix $\mathbf{T}:=\mathrm{blkdiag}(\Mmu^{1/2} ,\Meps^{1/2})$ 
allows to rewrite \eqref{eq:FIT_ode_1} as
\begin{align}
    \operatorname{d}_t \mathbf{u} &= \A \mathbf{u}+{\mathbf{g}}(t), & \mathbf{u}(t_0) &= \mathbf{u}_0
    \label{eq:FIT_ode_2}
\end{align}
in the new unknowns $\mathbf{u}^{\top}=[\Mmu^{1/2}\vh^{\top},\Meps^{1/2}\ve^{\top}]$
with the skew-symmetric stiffness matrix
\begin{align}
    \A=\begin{bmatrix}
        \mathbf{0} & -\Mmu^{-1/2}\curlfit\Meps^{-1/2}\\
        \Meps^{-1/2}\curldfit\Mmu^{-1/2} & \mathbf{0}
    \end{bmatrix}.
\end{align}
and right-hand-side $\mathbf{g}^{\top}=[\mathbf{0},\Meps^{-1/2}\fj_\text{s}^{\top}]$. 
One advantage of the transformed system \eqref{eq:FIT_ode_2} is that $\A$ is normal, 
i.e., $\A \A^{\top} = \A^{\top} \A$ and one shows straightforwardly that all eigenvalues 
are imaginary.

\subsection{Leapfrog}
\label{sec:Leapfrog}

For high-frequency electromagnetic initial value problems one typically employs 
the Leapfrog scheme (or equivalently St\"{o}rmer-Verlet) to solve the semi-discrete 
system \eqref{eq:FIT_ode_1}. If the initial condition 
\begin{align*}
	\ve^{(0)}=\ve_{0} \quad\text{and\quad} \vh^{(\frac{1}{2})}=\vh_{1/2}
\end{align*}
is given, the Leapfrog update equations read time step $m\in\{0,\ldots,n_{t}-1\}$
\begin{align*}
	\ve^{(m+1)}&=\ve^{(m)}+\Delta t \Meps^{-1}\left(\curldfit\vh^{(m+\frac{1}{2})}-\fj^{(m+\frac{1}{2})}\right),\\
	\vh^{(m+\frac{3}{2})}&=\vh^{(m+\frac{1}{2})}-\Delta t \Mmu^{-1} \curldfit\ve^{(m+1)}
\end{align*}
for the electric and magnetic voltages $\ve^{(m)}$, $\vh^{(m+\frac{1}{2})}$ at 
time points $t_m$ and $t_{m+\frac{1}{2}}$ with step size $\Delta t$. The scheme 
is (up to scaling) equivalent to Yee's Finite-Difference-Time-Domain scheme \cite{Yee}. 
Each equation applies the discretized curl operator and a few vector additions 
and scalar multiplications of complexity $\mathcal{O}(n_\text{dof})$. In the following discussion the focus will lie on the sparse-matrix-vector-multiplications (SMVP) as they are typically the most heavy operation, i.e., the cost of Leapfrog in terms of SMVP is given as 
\begin{align}
	\label{eq:cost_LF}
	C_\text{LF} = 2\cdot n_{t}
\end{align}
while additions and scalar multiplications are disregarded.

Leapfrog's time-stepping scheme is only conditionally stable and the maximal stable time step size 
$\Delta t$ is limited by the Courant-Friedrichs-Lewy (CFL) condition. A sharp 
bound is given by the largest absolute eigenvalue of $\A$, i.e., 
$\lambda_\text{max}=\| \A \|_2$ since $\A$ is normal. In practice, one estimates the value by 
\begin{align}
    \Delta t \leq \Delta t_\text{CFL} = \min\limits_j \left( \sqrt{\frac{\varepsilon_j\mu_j}{\frac{1}{\Delta x_j^2}+\frac{1}{\Delta y_j^2}+\frac{1}{\Delta z_j^2}}}\right), \label{eq:cfl}
\end{align}
where $j$ is the index of the grid cells and $\Delta x_j$ $\Delta y_j$ and 
$\Delta z_j$ their spatial dimensions.

One remarkable property of the Leapfrog scheme is energy preservation (or symplecticity). 
Let us define the discrete magnetic and electric energies as follows
\begin{equation}
	\begin{aligned}	
		\magenergy{m}{n}&:=\tr{\vh^{(m)}}\Mmu\vh^{(n)}\\
		\elecenergy{m}{n}&:=\tr{\ve^{(m)}}\Meps\ve^{(n)}.
	\end{aligned}	
\end{equation}
Using those energies and disregarding source currents, i.e. $\fj_\text{s}\equiv0$, one can show
\begin{equation}\label{eq:energy}
	\begin{aligned}	
		\elecenergy{m+1}{m+\frac{3}{2}}
		+
		\magenergy{m+1}{m+\frac{3}{2}}\\
		=
		\elecenergy{m}{m+\frac{1}{2}}
		+
		\magenergy{m}{m+\frac{1}{2}}
	\end{aligned}	
\end{equation}
where the electric and magnetic voltages must be consistently interpolated as 
$\ve^{(m+\frac{1}{2})}:=(\ve^{(m)}+\ve^{(m+1)})/2$ and 
$\vh^{(m)}:=(\vh^{(m-\frac{1}{2})}+\vh^{(m+\frac{1}{2})})/2$, see e.g. \cite{Bossavit}.     
\section{The ParaExp Algorithm}
\label{sec:paraexp}
In this section we develop based on \cite{gander2013paraexp} the ideas of the 
ParaExp method for the system of ODEs in the form \eqref{eq:FIT_ode_2}
\begin{align*}
	\operatorname{d}_t \mathbf{u} &=  \mathbf{\A} \mathbf{u}+\mathbf{g}(t) & \mathbf{u}(t_0)&= \mathbf{u}_0.
\end{align*}
Applying the method of variation of constants leads to the solution
\begin{equation}
	\mathbf{u}(t) = \exp{(t \A}) \mathbf{u}_0 + \int_0^t \exp{((t-\tau) 
	\A ) }\mathbf{g}(\tau) \mathrm{d} \tau
	\label{eq:solution-variation-of-constants}
\end{equation}
where $\exp{(t \A}) \mathbf{u}_0$ is the homogeneous solution due to initial 
conditions and the convolution product is the particular solution resulting from 
the presence of the source term $\mathbf{g}(t)$. The last term of 
\eqref{eq:solution-variation-of-constants} is more difficult to compute than the 
first one. However, thanks to the linearity of the problem and 
the superposition principle, $\mathbf{u}(t)$ can be written as 
$\mathbf{u}(t) = \mathbf{v}(t) + \mathbf{w}(t)$ where the particular solution 
$\mathbf{v}(t)$ is governed by
\begin{align}
	\operatorname{d}_t \mathbf{v} &= \A \mathbf{v}+\mathbf{g}(t) 
	& \mathbf{v}(t_0)&= \boldsymbol{0}
	\label{eq:solution-particular}
\intertext{and the homogeneous solution $\mathbf{w}(t)$ is governed by}
	\operatorname{d}_t \mathbf{w} &=\A \mathbf{w}+ \boldsymbol{0} 
	& \mathbf{w}(t_0)&= \mathbf{u}_0.
	\label{eq:solution-homogeneous}
\end{align}
The ParaExp method takes advantage of this decomposition. The time interval $\mathcal{I}=(0, T]$
is partitioned into sub-intervals $\mathcal{I}_j = (T_{\mathrm{j-1}}, T_{\mathrm{j}}]$
with $j = 1, 2, ..., p$,  $t_0 = T_0 < T_1 < T_2 < ... < T_p = t_{\mathrm{end}}$
and $p$ the number of CPUs. The following solutions are then computed on each CPU: 
a particular solution $\mathbf{v}_j:\mathcal{I}_j\to\real^{n_\text{dof}}$ governed by
\begin{align}
	\label{eq:solution-particular-paraexp}
	\operatorname{d}_t \mathbf{v}_j &= \A \mathbf{v}_j +\mathbf{g}(t) 
	&
	\mathbf{v}_j(T_{j-1})&= \boldsymbol{0}
\intertext{and a homogeneous solution $\mathbf{w}_j:\mathcal{I}_j\to\real^{n_\text{dof}}$ governed by}
	\label{eq:solution-homogeneous-paraexp}
	\operatorname{d}_t \mathbf{w}_j &= \A \mathbf{w}_j +\boldsymbol{0}
	&
	\mathbf{w}_j(T_{j-1})&= \mathbf{v}_{\mathrm{j-1}}(T_{j-1}).
\end{align}
Problems \eqref{eq:solution-particular-paraexp} can be solved in parallel by any 
time stepping method as only trivial initial conditions must be provided. In this paper 
Leapfrog is employed as discussed in Section~\ref{sec:Leapfrog}. The solutions 
for \eqref{eq:solution-homogeneous-paraexp} can be given analytically by 
\begin{align}
	\label{eq:matrix_expo}
	\mathbf{w}_j(t) = \exp{(t \A)} \mathbf{v}_{\mathrm{j-1}}(T_{j-1})
\end{align}
where the initial condition is the final solution $\mathbf{v}_{\mathrm{j-1}}(T_{j-1})$ of the previous interval. 
It is therefore highly recommended to compute $\mathbf{v}_{\mathrm{j-1}}(t)$ and 
$\mathbf{w}_{\mathrm{j}}(t)$ on the same CPU to avoid communicational costs.
Using the superposition principle, the total solution can be expanded as:
\begin{equation}
    \mathbf{u}(t) = \mathbf{v}_j(t) + \sum_{i = 1}^{j}\mathbf{w}_i(t) \, \, \textrm{ with } j \textrm{ s.t. } 
    t\in \mathcal{I}_j.
    \label{eq:solution-particular-para}
\end{equation}
Figure \ref{fig:timedecomposition} shows the time decomposition of IVP into 
particular solutions (solid lines) and homogeneous solutions (dashed lines) for 
a case with 3 CPUs. The two steps of the 
ParaExp method are described in the Algorithm shown in Figure~\ref{alg:ParaExp}.

\begin{algorithm}[t]
\SetKwInput{Input}{Input}
\SetKwInput{Output}{Output}
\SetKwProg{Parfor}{parfor}{ do}{end}
\LinesNumbered
\DontPrintSemicolon
\Input{system matrix $\mathbf{A}$, source term $\mathbf{g}(t)$, initial value $\mathbf{u}_0$, time interval $\mathcal{I}$, number of processors $p$}
\Output{solution $\mathbf{u}(t)$}
\Begin
{   partition $\mathcal{I}$ into intervals $\mathcal{I}_j$, $j=0,...,p$,\\
    \emph{\# begin the parallel loop (index $j$)}\\
    \Parfor{$(j \gets 1$ \KwTo $p)$}{
	$\mathbf{v}_j \gets$ solve $\mathbf{v}_j'(t)=\mathbf{A}\mathbf{v}_j(t)+\mathbf{g}(t)$,~$\mathbf{v}_j(T_{j-1})=0$, $t \in \mathcal{I}_j$ using a time stepper\\
	\eIf{$j\neq p$}{
        $\mathbf{w}_{j+1} \gets \exp(\mathbf{A}(t-T_{j}))\mathbf{v}_j(T_j)$ for all $t \in (T_j,T_p]$
        }{
        $\mathbf{w}_{1} \gets \exp(\mathbf{A}(t-T_{0}))\mathbf{u}_0$ for all $t \in (T_0,T_p]$
        }
    }
    \For{$(j \gets 1$ \KwTo $p)$}{
    $\mathbf{u}(t) \gets \mathbf{v}_j(t)+\sum\limits_{i=1}^j\mathbf{w}_i(t)$, for all $t\in \mathcal{I}_j$
    }
}
\caption{Pseudocode for the ParaExp Algorithm}
\label{alg:ParaExp}
\end{algorithm}
 
\subsection{Approximation of the matrix exponential}
\label{sec:paraexp_matrix_exponential}
A critical point of the method is the efficient computation of equation \eqref{eq:matrix_expo}
by the matrix exponential. A 
straight forward evaluation of the exponential followed by the multiplication with 
the vector of initial conditions is computationally very costly, especially for 
large matrices such as the matrices obtained by spatial discretization of the wave problem. 
Instead, efficient approximations of the action of the matrix exponential on initial 
condition vectors are used.
Examples of such methods are the Krylov subspace based methods as used in \cite{gander2013paraexp}, 
Higham's function \cite{alhi11} and Leja's method \cite{caliari2016leja}.

The Krylov subspace based methods (see \cite{gander2013paraexp}) require the evaluation 
of Krylov subspaces $\mathcal{K}^l(\boldsymbol{S}, \bb):=\left\{ \bb, \boldsymbol{S}\bb, \cdots, \boldsymbol{S}^{l-1}\bb \right\}$. 
These subspaces involve the evaluation of the matrix 
\begin{equation}
    \boldsymbol{S} := \left( \boldsymbol{I} - \A/\sigma \right)^{-1} \A \in \complex^{2n_\text{dof} \times 2n_\text{dof}}
    \label{eq:krylov_subspace}
\end{equation}
with $\sigma\in\complex$ and its multiplication with the vector $\bb$ which is for 
example given by some solution $\mathbf{v}_{j}(t)$. If $\sigma \neq \infty$, one 
may approximate the action of the exponential with a rather small $n$ but the 
computational costs of solving the large linear systems in \eqref{eq:krylov_subspace} 
become prohibitive. The choice $\sigma = \infty$ leads to 
$\mathcal{K}^n(\boldsymbol{S}, \bb) = \mathcal{K}^n(\boldsymbol{A}, \bb)$ and 
avoids matrix inversions but typically requires a rather large Krylov subspace 
in practice. In either case, the computational cost associated with the Krylov 
subspace based methods have been rather large such that we focus on Higham's 
and Leja's method in the following. 

Both methods use two main ingredients.
The first ingredient is the \emph{scaling} of the matrix exponential:
\begin{align*}
    \exp \left( t \A \right) \bb &= \left( \exp \left(t \A/s \right) \right)^s \bb \\
    &\approx \left( P \left(t \A/s \right) \right)^s \bb 
		=: \exp \left( t \left(\A + \Delta \A \right) \right) \bb
\end{align*}
which reduces the spectrum of the scaled matrix $\exp \left(t \A/s \right)$ around 
the origin thus allowing its efficient approximation by polynomial interpolations
$P$ such as Taylor's expansion.
The second ingredient is the use of a \emph{recurrence equation} that involves 
SMVPs. 
In the case Higham's function, the recurrence equation reads:
\begin{align*}
    \bb_{i+1} = T_m\left( t\A/s \right) \bb_i, i = 0, \cdots ,s-1, \bb_0 = \bb.
\end{align*}
where $T_m$ is Taylor's truncated polynomial of order $m$, i.e., 
\begin{align*}
    T_m \left( t\A/s \right) = \sum \limits_{j=0}^{m}\frac{\left(t \A/s \right)^j}{j!}.
\end{align*}
Leja's method uses interpolation which is a Newton-Cotes interpolation polynomial
$L_{m, c}$ defined on the set of Leja's points \cite{caliari2016leja}. Similarly to
Chebyshev's approach, the points are chosen such that the condition number of the
polynomial remains small when the polynomial order $m$ is increasing. In the case 
of the wave equation with imaginary eigenvalues, the interpolation is defined on an 
interval in the complex plane. The resulting recurrence equation reads
\begin{align*}
    \bb_{i+1} = L_{m, c} \left(t \A/s \right) \bb_i, i=0,\cdots, s-1, \bb_0 = \bb. 
\end{align*}
Leja's interpolation $L_{m, c}$ reduces to the Taylor series for $c = 0$ but it 
performs better than Higham's functions for normal matrices with a value of $c \neq 0$.
In both cases, the approximate solution is derived as $\bu = \bb_{s}$. 
The parameters $m, s$ (and $c$ in the case of Leja's method) are chosen so as to 
minimize the computational cost given by the number of SMVPs
$$
C_\text{Leja} = n_\text{Leja} ,
$$
with the condition $\| \Delta \A \| \leq \varepsilon_A \| \A \|$ where $\varepsilon_A$ 
is a prescribed tolerance. 
{This cost is dominated by the approximation of the matrix exponential with
$n_\text{Leja} \approx s \, m$. The approximation involves $s$ multiplications 
of the Taylor polynomial $P \left(t \A/s \right)$ of order $m$. An additional cost 
results from the evaluation of the optimal parameters $m, s$ and $c$ involving the 
computation of the norm of the matrix A (see \cite{alhi11}).
}
Both algorithms also use additional preprocessing steps (shifting and balancing of the matrix)
and an early termination of the iteration in the polynomial interpolations.

\subsection{Reconstruction of voltages on staggered grids}
\label{sec:paraexp_staggered}

When solving \eqref{eq:solution-particular-paraexp} with Leapfrog, electric and 
magnetic voltages are allocated on staggered time grids and it has been shown in 
Section~\ref{sec:Leapfrog} that this is crucial for energy conservation. {Therefore 
we propose to modify the reconstruction \eqref{eq:solution-particular-para} in 
the case of Leapfrog in the obvious manner for the total electric grid voltages as
\begin{align}
	\label{eq:averaging1}
	\ve^{(m+1)}&=\ve^{(m+1)}_j +[\mathbf{0},\Meps^{-1/2}] \sum_{i = 1}^{j}\mathbf{w}_i(t^{m+1})
\intertext{for $t^{m+1}\in\mathcal{I}_j$ and the total magnetic grid voltages}
	\label{eq:averaging2}
	\vh^{(m+\frac{1}{2})}&=\vh^{(m+\frac{1}{2})}_j +[\Mmu^{-1/2},\mathbf{0}] \sum_{i = 1}^{j}\mathbf{w}_i(t^{m+\frac{1}{2}})
\end{align}
for $t^{m+\frac{1}{2}}\in\mathcal{I}_j$ where $\ve^{(m+1)}_j$ and $\vh^{(m+\frac{1}{2})}_j$ are the solutions of 
\eqref{eq:solution-particular-paraexp} at time $t^{m+1}$ and $t^{m+\frac{1}{2}}$ 
using Leapfrog and the matrices $[\mathbf{0},\Meps^{-1/2}]$ and $[\Mmu^{-1/2},\mathbf{0}]$ 
are used to extract and transform the respective components from the solution $\mathbf{w}$ of \eqref{eq:solution-homogeneous-paraexp}. Leapfrog is initialized with $\ve_0=\mathbf{0}$ and $\vh_\frac{1}{2}=\mathbf{0}$ for all time intervals $\mathcal{I}_j$.}
\subsection{Discussion of the numerical costs}
\label{sec:paraexp_costs}
The computational costs of Paraexp can be split into three categories: (i) execution of Leapfrog, (ii) propagation of initial values by the matrix exponential and (iii) two transformations \eqref{eq:FIT_ode_2} for initial and end values at each interval. The effective number of SMVPs, i.e., disregarding operations carried out in parallel, is given by
\begin{equation}
	\label{eq:tot_cost}
	C_\text{proc} = \frac{2}{p}n_{t} + n_\text{Leja} + 2,
\end{equation}
if the same number of time steps is performed on each processor. When increasing the number or processors $p$, the costs of Leapfrog can be disregarded and the only remaining costs are SMVPs due to Leja for the longest time interval, i.e., for propagating the initial value from $t_0$ to $t_\mathrm{end}$.
\subsection{Error analysis}
\label{sec:paraexp_error_analysis}
Neglecting the round-off errors, the contribution to the main error of ParaExp method 
from the $j^{\mathrm{th}}$ CPU has two main contributions. 

The first contribution is the truncation error resulting 
from the Leapfrog scheme used for solving the non-homogeneous problem on the sub-interval 
$\mathcal{I}_j=(T_{j-1},T_{j}]$. There is no error in the initial value since we start 
from the trivial initial condition. The resulting numerical solution $\hat{\bv}_{j}$ 
is given by
\begin{equation}
    \hat{\bv}_{j} = \bv_{j} + \Delta \bv_{j},
\end{equation} 
where $\bv_{j}$ is the exact solution of the problem and $\Delta \bv_{j} \sim \mathcal{O} \left( \Delta t \right)^2$ 
is the truncation error. This solution is then used as an initial condition for 
the homogeneous problem. 

The second contribution is the approximation error of the matrix exponential on 
the interval $(T_{j}, T]$. It can be analyzed using the backward error analysis 
as introduced by Higham, see Section \ref{sec:paraexp_matrix_exponential}:
\begin{equation}
    \hat{\bw}_{j} = \bw_{j} + \Delta \bw_{j} = \exp\left(t \left( \A + \Delta \A \right) \right) \hat{\mathbf{v}}_{\mathrm{j}}(T_{j})
\end{equation}
with $\hat{\bw}_{j}$ the numerical solution obtained by approximating 
the exponential of the matrix, $\bw_{j}$ the exact homogeneous 
solution and $\Delta \bw_{j}$ the approximation error. 
The contribution of the $j^{\mathrm{th}}$ CPU to the total numerical solution is 
therefore given by
\begin{align*}
    \hat{\bu}_{j}  
    &= \bu_{j} + \Delta \bu_{j}
    = \exp\left(t \A \right) \bv_{j} + \Delta \bu_{j}
    \\
    &=\exp\left(t \left( \A + \Delta \A \right) \right) \left(\bv_{j} + \Delta \bv_{j} \right),
\end{align*}
where  $\bu_{j}=\bv_{j}+\bw_{j}$ is the exact solution.
The analysis developed in Section 4 of \cite{alhi11} can be used to quantify these
errors and to adjust the relative tolerance of Leja's method so that both errors 
are of the same magnitude. Applying Lemma 4.2 from \cite{alhi11} together with formula
(4.5) to the action of the matrix exponential with $\| \Delta \A \| \leq \varepsilon_A \| \A \|$ 
where $\varepsilon_A$ is a prescribed relative tolerance and assuming that 
$\| \Delta \bv_{j} \| \leq \varepsilon_B \| \bw_{j} \|$ with 
$\varepsilon_B = \beta {\Delta t}^2$, the following result can be derived:
\begin{equation*}
    \frac{\| \Delta \bu_{j} \|_2}{ \| \bu_{j} \|_2} 
    \leq 
    \frac{\| e^{t\A}\|_F \| \bv_{j} \|_2}{ \| \bu_{j} \|_2} 
    \left( \frac{\varepsilon_B}{\varepsilon_A} + \kappa_{\mathrm{exp}}\left(\A \right) \right)
\end{equation*}
where $\kappa_{\mathrm{exp}} \left( \A \right)$ is the condition 
number of the matrix exponential, i.e., $\kappa_{\mathrm{exp}}(\A) = \| \A\|_2$ for 
normal matrices. For a fixed $\varepsilon_B$, 
both terms in the brackets become equal if
\begin{equation}
    \frac{\varepsilon_B}{\varepsilon_{A}^{\ast}} = \kappa_{\mathrm{exp}}(\A) = \| \A \|_2
    \label{eq:tolerance_leja_1}
\end{equation}
leading to the optimal value of the tolerance for the matrix exponential given by
\begin{equation}
    \varepsilon_{A}^{\ast} = \frac{\beta {\Delta t}^2}{\kappa_{\mathrm{exp}}(\A)} = \frac{\beta {\Delta t}^2}{\| \A \|_2}.
    \label{eq:tolerance_leja_2}
\end{equation}

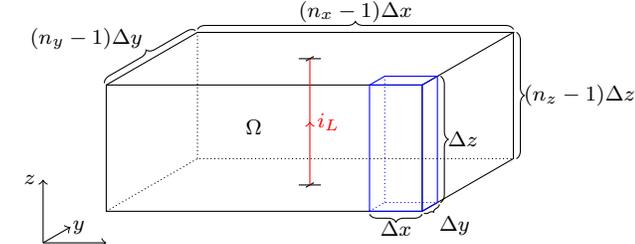
\begin{figure}[t]
    \centering
    \usetikzlibrary{decorations.markings}
\usetikzlibrary{decorations.pathreplacing}
\usetikzlibrary{arrows}
\begin{footnotesize}
\tikzset{
  c/.style={every coordinate/.try}
}
\begin{tikzpicture}[scale=1.4,nodewitharrow/.style 2 args={decoration={markings,mark=at position {#1} with { \arrow{>},\node[transform shape,above] {#2};}},postaction={decorate}}]
\def\h{1.2};
\def\w{3};
\def\dh{1.2};
\def\dw{0.5};
\def\a{30};
\def\denom{3};

	\draw (0,0)coordinate[name=p1]{} rectangle (\w,\h)coordinate[name=p4]{};
	\draw (0,\h)coordinate[name=p3]{}--++(\a:\w/\denom)coordinate[name=p7]{}--++(0:\w)coordinate[name=p8]{}--++(-180+\a:\w/\denom);
	\draw (\w,0)coordinate[name=p2]{}--++(\a:\w/\denom) coordinate[name=p6]{} --++(90:\h) coordinate[name=p8]{};
	\draw [densely dotted] (p1)--++(\a:\w/\denom)coordinate[name=p5]{}--(p6);
	\draw [densely dotted] (p5)--(p7);
	\path (\w/2,0)--++(\a:\w/\denom/2)coordinate[name=mu]{}--++(0,\h)coordinate[name=mo]{};

    \draw[thin] (mu)--++(\a:1.25pt); \draw (mu)--++(-180+\a:1.25pt); \draw (mu)--++(3pt,0); \draw (mu)--++(-3pt,0); 
	\draw [thin] (mo)--++(\a:1.25pt); \draw (mo)--++(-180+\a:1.25pt); \draw (mo)--++(3pt,0); \draw (mo)--++(-3pt,0); 
	\draw [red,nodewitharrow={0.5}{}] (mu)--(mo) node [pos=0.5,right]{$i_L$};

\begin{scope}[shift={(\w-\dw,0)}]
\draw [blue](0,0) coordinate[name=dp1]{} rectangle (\dw,\h);
\draw[blue] (0,\dh)  coordinate[name=dp3]{}--++(\a:\dw/\denom) coordinate[name=dp7]{}--++(\dw,0) coordinate[name=dp8]{}--(\dw,\dh) coordinate[name=dp4]{};
\draw [blue, densely dotted] (0,0)--++(\a:\dw/\denom) coordinate[name=dp5]{}--++(\dw,0) coordinate[name=dp6]{};
\draw[blue, densely dotted] (dp5)--++(0,\dh);
\draw[blue] (\dw,0)coordinate[name=dp2]{}--(dp6)--++(0,\dh);
\end{scope}
\begin{scope}[scale=0.2,xshift=-3cm,yshift=-1.5cm]
		\draw [->] (0,0)--(3,0) node[below] {$x$};
		\draw [->] (0,0)--(0,3)node[left] {$z$};
		\draw [scale=1.5,->] (0,0)--(\a:1)node[right,xshift=-2pt] {$y$};
\end{scope}
	\node (1) at (1.4,0.8) {$\Omega$};

\begin{scope}[every coordinate/.style={shift={(0pt,0.8pt)}}]
 \draw  [decorate,decoration={brace,amplitude=3pt}]([c]p7)--([c]p8) node [yshift=0pt,midway, above] {$(n_x-1)\Delta x$};
\end{scope}
\begin{scope}[every coordinate/.style={shift={(-0.5pt,0.5pt)}}]
 \draw  [decorate,decoration={brace,amplitude=3pt}]([c]p3)--([c]p7) node [yshift=0pt,midway, above left] {$(n_y-1)\Delta y$};
\end{scope}
\begin{scope}[every coordinate/.style={shift={(0.8pt,0pt)}}]
 \draw  [decorate,decoration={brace,amplitude=3pt}]([c]p8)--([c]p6) node [yshift=0pt,midway,right] {$(n_z-1)\Delta z$};
\end{scope}
\begin{scope}[every coordinate/.style={shift={(0pt,-0.8pt)}}]
 \draw  [decorate,decoration={brace,amplitude=3pt}]([c]dp2)--([c]dp1) node [yshift=0pt,midway, below] {$\Delta x$};
\end{scope}
\begin{scope}[every coordinate/.style={shift={(0.5pt,-0.5pt)}}]
 \draw  [decorate,decoration={brace,amplitude=1.5pt}]([c]dp6)--([c]dp2) node [yshift=0pt,midway, below right] {$\Delta y$};
\end{scope}
\begin{scope}[every coordinate/.style={shift={(0.8pt,0pt)}}]
 \draw  [decorate,decoration={brace,amplitude=3pt}]([c]dp8)--([c]dp6) node [yshift=0pt,midway,right] {$\Delta z$};
\end{scope}
	
\end{tikzpicture}
\end{footnotesize}
     \vspace*{-2em}
    \caption{Domain $\Omega$ of the wave problem with a hexahedral mesh}
    \label{fig:waveProbDef}
\end{figure}     
\section{Numerical tests}
\label{sec:applications}

We consider two numerical tests for the validation: a two-dimensional cylindrical 
wave excited by a line current and a spiral inductor discretized by \cstmws based on the design proposed in \cite{Becks_1992aa,CST_2016aa}. 

\subsection{Cylindrical two-dimensional wave}
\label{subsec:applications}

The two-dimensional cylindrical wave problem is depicted in Figure \ref{fig:waveProbDef}. 
The excitation is a line current in $z$-direction in the center of the domain $\Omega$. The discretization is obtained by FIT as explained in Section \ref{sec:maxwell-equations}. 
A PEC boundary is assumed on the whole boundary $\partial \Omega$. The dimensions of the domain $\Omega$ are $L_x=L_y=\SI{20}{m}$ and $L_y=\SI{1}{m}$. For the discretization we use 
$n_x = n_y = 41, \, n_z = 2$ for Leapfrog and Leja's method 
and $n_x = n_y = 121, \, n_z = 2$ for the reference numerical solution. 
This corresponds to $n_\text{dof}=\num{20172}$ and $n_\text{dof}=\num{175692}$ degrees of freedom, respectively.

\begin{figure}[t]
    \centering
    \begin{footnotesize}
\begin{tikzpicture}[scale=1]
  \begin{axis}[xlabel={$x$ \si{[m]}},ylabel={$y$ [m]},zlabel={$\ve_z$},shader=faceted, yscale=0.6,scale=0.9,ytick={0,10,20,30}]
    \addplot3[surf,mesh/ordering=x varies,mesh/rows=31,colormap/jet,shader=faceted interp] table [x index=0, y index=1, z index=2, col sep=comma] {data/wavePlot2.csv};
  \end{axis}
\end{tikzpicture}
\end{footnotesize}
     \vspace*{-1em}
    \caption{The $\ve_z$ component of the wave at $t=\SI{4.8e-8}{\s}$}
    \label{fig:wave}
\end{figure}

The domain is filled with vacuum. The line current is a Gaussian function given by
\begin{align}
    i_\text{L}(t) = i_\text{max}\;\exp\left({-4\left(\displaystyle \frac{t-\sigma_t}{\sigma_t}\right)^2}\right) \label{eq:line_current}
\end{align}
with $i_\text{max}=\SI{1}{\A}$ and $\sigma_t=\SI{2e-8}{\s}$.
 
The differential equation of this problem is given by \eqref{eq:FIT_ode_1} and 
\eqref{eq:fit_matrices} with $\bar{\mathbf{g}}(t) = -[\boldsymbol{0},\fj]^{\top}$, 
with $\fj$ being the discretized line current \eqref{eq:line_current}.
We consider the a transformed ODE \eqref{eq:FIT_ode_2} with $\A$ 
being normal. The $\ve_z$ component of the calculated wave can be seen in Figure \ref{fig:wave}  
at time $t=\SI{4.4e-8}{\s}$.

Numerical experiments show that Leja's method outperforms Higham's approach and 
the Krylov subspace methods by a factor of $2$ and $10$, respectively. 
This is in agreement with the literature where Leja has been observed to perform 
best for normal matrices \cite{caliari2016leja}. Therefore, we will only present 
results of Leja's method from now on.

Leapfrog and ParaExp cause two different kinds of error: the truncation error for 
Leapfrog is related to the time step whereas the error of ParaExp can be quantified 
using the backward error analysis as discussed above. Therefore it is crucial to 
choose the parameters of Leja's algorithm to make the comparison between Leapfrog 
and ParaExp as fair as possible. 

For a given spatial mesh, the time step for Leapfrog is chosen according to 
the CFL criterion \eqref{eq:cfl}. The tolerance $\varepsilon_A$ of Leja's algorithm 
has been increased such that the error $\Delta\bv$ due to Leapfrog still dominates 
over the error due to the matrix exponential $\Delta\bw$. Results are shown in Figure 
\ref{fig:cost_LF_Leja_uniform_mesh}. 
In this figure, the costs of the Leapfrog scheme and Leja
are compared on $\mathcal{I}=(0, T]$ with $T =\SI{2e-7}{s}$ motivated by the  reasoning given in \eqref{eq:tot_cost}: if $p \gg 1$ we neglect the cost of the Leapfrog scheme used for solving 
the non-homogeneous problem \eqref{eq:solution-particular-paraexp} and only consider the cost of the propagation of the matrix exponential on $\mathcal{I}$.
As expected, the number of SMVPs of the Leapfrog algorithm is linearly proportional 
to the number of time steps $n_t$, see \eqref{eq:cost_LF}. 
The increase of the number of time steps $n_t$ does not correspond to an increase 
of cost for Leja method as the matrix of the system and the time interval remain unchanged. 
A slight increase may however result from the evaluation of intermediate interpolations.

\begin{figure}[t]
    \begin{tikzpicture}
        \begin{axis}[zlabel={SMVPs},ylabel={$n_t$},xlabel={$n_x$}, width=0.9\columnwidth,xtick scale label code/.code={},
        zmin=0,zmax=8000,shader=faceted interp,point meta max=8000, point meta min=0,ytick={200,500,800},grid=major]
            \addplot3[surf,mesh/ordering=x varies,mesh/rows=19] table [x index=0, y index=1, z index=2, col sep=comma] 
            {data/nxntnSMVP_LF.csv};
        \end{axis}
    \end{tikzpicture}
    \begin{tikzpicture}
        \begin{axis}[zlabel={SMVPs},ylabel={$n_t$},xlabel={$n_x$}, width=0.9\columnwidth,xtick scale label code/.code={},
        zmin=0,zmax=8000,shader=faceted interp,point meta max=8000, point meta min=0, ytick={200,500,800},grid=major]
            \addplot3[surf,mesh/ordering=x varies,mesh/rows=19] table [x index=0, y index=1, z index=4, col sep=comma] 
            {data/nxntnSMVP_LF_Hig_Leja1_Leja2_Leja4.csv};
        \end{axis}
    \end{tikzpicture}
    \caption{Computational cost of Leapfrog and Leja's method for a uniform mesh.
            Top: cost for Leapfrog. Bottom: cost for Leja's method.
            The tolerance of Leja has been chosen as $\num{1e-2}$}
    \label{fig:cost_LF_Leja_uniform_mesh}
\end{figure}
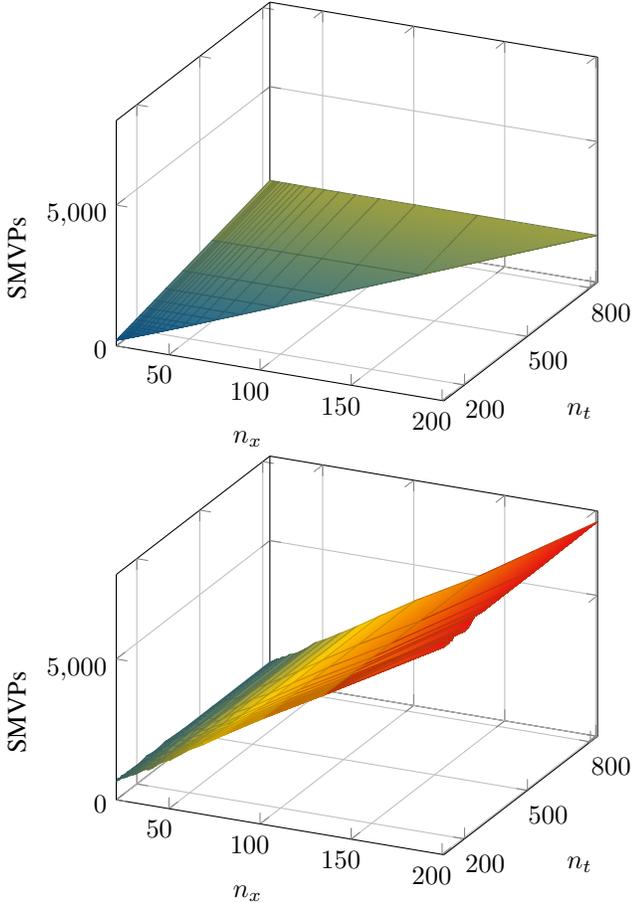 

The increase of $n_x$ which corresponds to the refinement of the spatial grid does 
not change the number of SMVPs for Leapfrog if ${\Delta t}$ remains below 
${\Delta t}_\text{CFL}$. Otherwise it increases linearly which is in good agreement 
with the estimate \eqref{eq:cfl}. For Leja's method, an even stronger linear increase 
of the number of SMVPs is observed also corresponding to the increase of the largest 
eigenvalue \eqref{eq:tolerance_leja_1}. In both cases the cost of each SMVP also 
increases due to the growing dimension $n_\text{dof}$ of the sparse matrices.

The Leapfrog scheme is known to be very competitive when using homogeneous grids, 
however small elements may quickly deteriorate the efficiency. Therefore a second 
numerical experiment investigates the efficiency in the cases of increasing inhomogeneity. 
The mesh of Figure \ref{fig:waveProbDef} was kept except for one element whose 
dimensions have been modified so that the ratio $k$ between the size of the biggest 
element over the size of the smallest element of the mesh lies in the interval 
$[1 , 20]$. The ratio of the computational cost
\begin{equation}
    R(k) = \frac{C_\text{Leja}(k)}{C_\text{LF}(k)}
    \label{eq:ration-costs}
\end{equation}
is given as the quotient of the number of SMVPs for performing Leja and Leapfrog 
on the whole time interval $\mathcal{I}$ in dependency of the non-uniformity of the mesh. 
Figure \ref{fig:cost_LF_Leja_nonuniform_mesh} shows a better performance for ParaExp over Leapfrog for large values of $k$. 
This suggests a better performance of ParaExp for highly non-uniform grids, although 
also Leja is depending via $\kappa_{\mathrm{exp}}$ on the eigenvalues of the operator $\mathbf{A}$.

\begin{figure}[t]
    \begin{tikzpicture}
	\begin{axis}[width=0.5\textwidth,height=6cm,xmin=1,xmax=20,ymin=0,ymax=4.25,xlabel={$k$},ylabel={$R$},y label style={at={(axis description cs:0.07,0.5)}}]
	   \addplot+[color=blue,mark=*,mark size=1.5,mark options=solid] table [x=a, y=b, col sep=comma] {data/ratioFLOPs.csv};
	   \addplot+[color=red,,mark=none,mark size=0.7,mark options=solid] table [x=a, y=c, col sep=comma] {data/ratioFLOPs.csv};
	   \legend{Ratio of computational costs}
	\end{axis}
    \end{tikzpicture}
    \vspace{-2em}
    \caption{Ratio between computational costs of Leapfrog and Leja's method for 
    a non-uniform mesh and comparable accuracies as defined in \eqref{eq:ration-costs}.}
    \label{fig:cost_LF_Leja_nonuniform_mesh}
\end{figure}
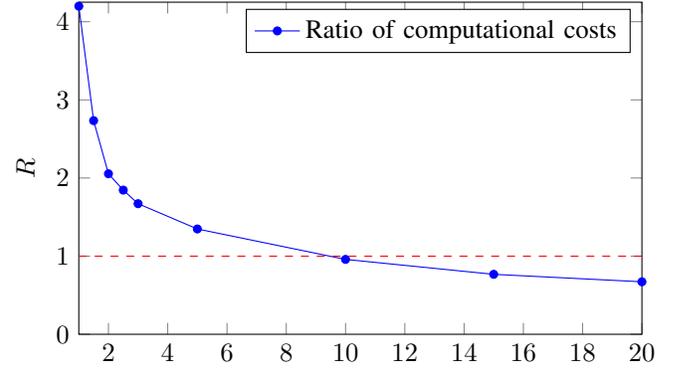
In a third experiment the energy conservation is numerically analyzed. ParaExp 
itself is as energy conserving as the methods used for time integration and propagation. 
While this is well understood for Leapfrog as discussed in Section \ref{sec:Leapfrog}, 
it is not clear for the Leja method. Figure \ref{fig:symplecticity} shows the energy for Leapfrog computed according to \eqref{eq:energy} and for Paraexp using
\begin{equation*}
	E(t):=
	\bigl\langle\vh(t),\vh(t)\bigr\rangle_\mu
	+
	\bigl\langle\ve(t),\ve(t)\bigr\rangle_\varepsilon
\end{equation*}
with the averaging from \eqref{eq:averaging1}.
It can be seen that the energy remains constant once the excitation is vanishing, i.e., for $t \in \mathcal{I}$ such 
that $i_L(t) = 0$. The electromagnetic energy present in the system remains constant 
and almost independent of the spatial refinement. 

In order to analyze the accuracy of the ParaExp method the frequency spectrum of the electric field obtained by ParaExp and Leapfrog are compared in Figure \ref{fig:frequencyE}. It can be observed that ParaExp adds high frequency noise to the solution. {A side effect of those high frequencies is a potentially unphysical increase of energy in the domain as observed in Fig.~\ref{fig:symplecticity}. This can be avoided by reducing the time step size $\Delta t$ of the used Leapfrog algorithm. In this case using $\Delta t=\frac{\Delta t_\text{CFL}}{5}$ leads to a solution without high frequency noise.}
\begin{figure}
    \centering
    \pgfplotsset{cycle list={
tud1d ,solid\\
tud9d,dashed\\
tud5b ,densely dashed \\
violet,solid\\}}
\begin{tikzpicture}
\begin{semilogyaxis}[width=0.9*\columnwidth,ylabel={$\ve$ [V]},xlabel={frequency [Hz]},xmin=0,xmax=2.5e8,legend entries = {reference (Leapfrog),Leapfrog,ParaExp (Leja)}]
		\addplot table[x index=0,y index=1,col sep=comma] 
			{data/FD_e_freqRefLFLeja.csv};		
		\addplot table[x index=0,y index=2,col sep=comma] 
			{data/FD_e_freqRefLFLeja.csv};
		\addplot table[x index=0,y index=3,col sep=comma] 
			{data/FD_e_freqRefLFLeja.csv};
\end{semilogyaxis}
\end{tikzpicture}
     \caption{Frequency spectrum of $\ve$ for Leapfrog and ParaExp (6 parallel threads, Leja).}
    \label{fig:frequencyE}
\end{figure}
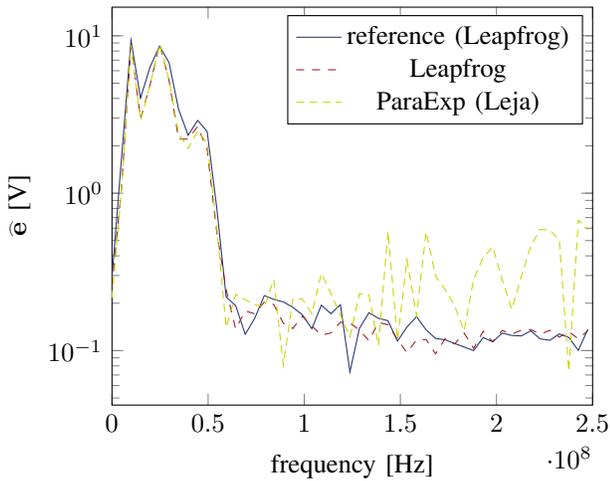

\subsection{Spiral inductor\label{sec:cst}}
The second test case is a spiral inductor model created with coplanar lines located 
on a substrate layer with an air bridge. The model is based on the design proposed 
in \cite{Becks_1992aa} and the corresponding example from the CST tutorial on transient 
analysis which advocates the usage of 3D field simulation instead of circuit models, 
cf. \cite{CST_2016aa}. The geometry of the problem is illustrated in Figure \ref{fig:SpiralInductor_geometry} 
with the dimensions: $\SI{7e-4}{m} \times \SI{4.75e-4}{m} \times \SI{2.5e-5}{m}$. 
The structure is discretized using a mesh with $\num{406493}$ mesh cells. The 
largest cell has the dimensions $\SI{8.9e-6}{m} \times \SI{8.9e-6}{m} \times \SI{8.8e-6}{m}$ 
and the smallest cell has the dimensions $\SI{2e-6}{m} \times \SI{2e-6}{m} \times \SI{1.5e-6}{m}$. 
The conductor is modeled by PEC, the substrate is given by a relative permittivity of $\epsilon_\text{r}=12$.
The domain is discretized using $n_\text{dof}=1,283,040$ degrees of freedom for 
$\ve$ and $\vh$ respectively and excited by a sine wave at $\SI{50}{GHz}$. The 
propagation of the results from $t_0$ to $t_\text{end}$ by Leapfrog requires $\num{21654}$ SMVPs 
while Leja needs $\num{34864}$ SMVPs for the same interval, cf. Figure \ref{fig:SpiralInductor_results}. 
In other words, classical time-stepping by Leapfrog is approximately $1.6$ times faster
than the evaluation of the matrix exponential for this example and ParaExp does not pay off. 
Motivated by the results of Figure~\ref{fig:cost_LF_Leja_nonuniform_mesh}, a carefully chosen 
example with small geometric details may change the situation in favor of ParaExp but for 
a general problem Leapfrog remains the appropriate choice.
\begin{figure}[t]
	\hspace*{-1em}
	\begin{tikzpicture}
		\begin{semilogyaxis}[width=0.48\textwidth,legend style={at={(axis cs:1.97e-7,2e-9)},anchor=north east},height = 6cm,xmin = 0,xmax = 2e-7,ymin = 0,xlabel = {Time (s)}, ylabel={Energy (J)}]
			\addplot+[color=red, dashed, thick, mark = *, mark size = 0.7, mark options = {solid}, smooth] table [x = a, y = d, col sep = comma] {data/energy_t_Leja_LFRef_LFLF.csv};
			\addplot+[color=black, dashed, thick, mark = *, mark size = 0.7, mark options={solid}, smooth] table [x = a, y = b, col sep = comma]  {data/energy_t_Leja_LFRef_LFLF.csv};
			\legend{Leapfrog, ParaExp}
		\end{semilogyaxis}
	\end{tikzpicture}
  \vspace{-2em}
	\caption{Electromagnetic energy obtained for Leapfrog and ParaExp using Leja.}
	\label{fig:symplecticity}
\end{figure}
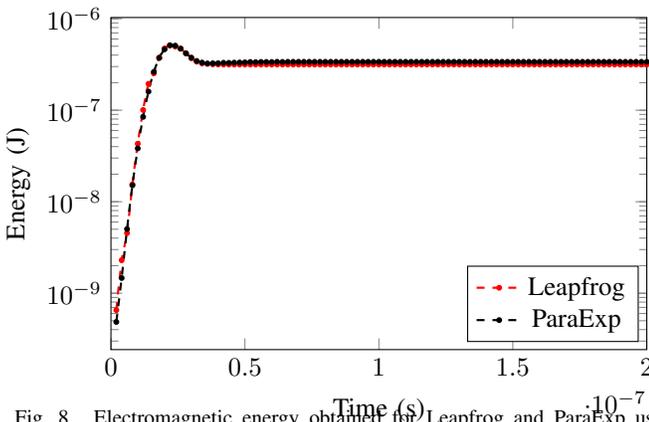
     
\section{Conclusions\label{sec:conclusions}}
In this paper, the ParaExp method was used for parallelization of time domain simulation 
of the electromagnetic wave problem and its performances were compared to the 
performance of the Leapfrog scheme.
The efficiency of ParaExp heavily depends on the approximation of the action of 
the matrix exponential to the vector of initial conditions and this approximation 
is the bottleneck of the method. Two methods were investigated for the approximation 
of the matrix exponential: Higham's function and Leja's method. In our applications, 
Leja's method was more efficient because the problem can be reformulated in terms of a normal system matrix. 
Numerical tests have shown that Leapfrog performs better than ParaExp with Leja for problems involving uniform meshes, but ParaExp can become more efficient than Leapfrog for problems involving highly non-uniform meshes.
In contrast to Leapfrog the energy preservation of Paraexp does not only depend on the time stepper used but also on the approximation of the matrix exponential.     \begin{figure}[t]
	\centering
	\includegraphics[trim={0 0 0 1cm},width=0.95\columnwidth]{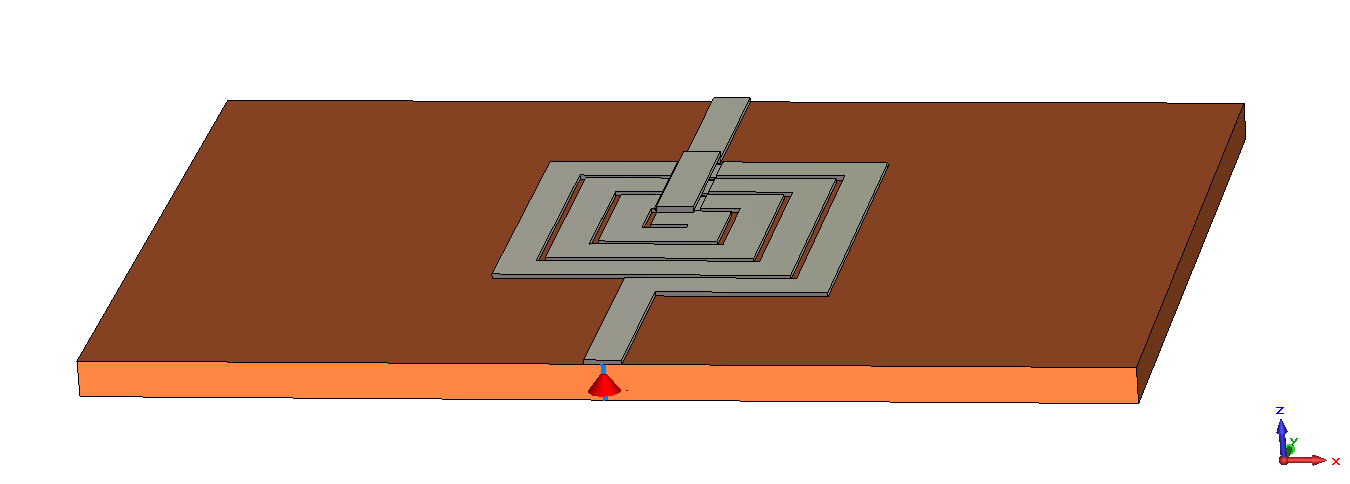}
	\vspace{-1em}
	\caption{\cstmws model of a spiral inductor based on \cite{Becks_1992aa}}
	\label{fig:SpiralInductor_geometry}
\end{figure}
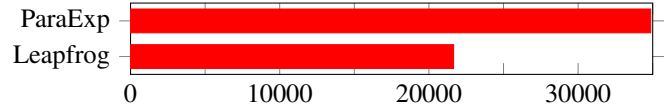
\begin{figure}[t]
  \hspace*{-1.5em}\begin{tikzpicture}
      \begin{axis}[xbar, width=0.47\textwidth, bar width=2ex, y=3ex, enlarge y limits={abs=0.5},yticklabels={Leapfrog,ParaExp},ytick={1,2},xtick={0,5000,10000,15000,20000,25000,30000,35000},xticklabels={0,,10000,,20000,,30000,},scaled x ticks=false,xmin=0,,xmax=35000]
          \addplot[red,fill=red] coordinates {(21654, 1) (34864, 2)};
      \end{axis}
      \label{fig:SpiralInductor_results}
  \end{tikzpicture}
	\vspace{-0.7em}
	\caption{Number of SMVPs of Leja and Leapfrog to propagate the initial value from $t_0$ to $t_\text{end}$}
\end{figure}

\section*{Acknowledgments}
The authors would like to thank Timo Euler, CST AG for the fruitful discussions 
on time domain simulations and Prof. Marco Caliari from the University of Verona 
for the discussion on the numerical implementation of Leja's method. 

This work was supported by the German Funding Agency (DFG) by the grant 
`Parallel and Explicit Methods for the Eddy Current Problem' (SCHO-1562/1-1), 
the 'Excellence Initiative' of the German Federal and State Governments and the 
Graduate School CE at Technische Universit\"at Darmstadt.
}


\end{document}